\documentclass[10pt,onecolumn,article]{IEEEtran}
\usepackage{times}
\usepackage{amsbsy}
\usepackage{latexsym}
\usepackage{amssymb}
\usepackage{amsmath,amssymb,amsfonts,graphicx,epsfig}
\usepackage{times,amsbsy,latexsym,bm,color}
\usepackage[ansinew]{inputenc}

\title{Technical Report on Two-Step Knowledge-Aided Iterative ESPRIT Algorithm \vspace{-0.01em}}
\author{Silvio F. B. Pinto $^1$ and Rodrigo C. de Lamare $^{1,2}$ \\
Center for Telecommunications Studies (CETUC) \\ $^1$ Pontifical
Catholic
University of Rio de Janeiro, RJ, Brazil.\\
$^2$ Department of Electronics, University of York, UK \\
Emails: silviof@cetuc.puc-rio.br, delamare@cetuc.puc-rio.br
\vspace{-0.01em}}
\begin{document}
\maketitle
\begin{abstract}
In this work, we  propose a subspace-based algorithm for
direction-of-arrival (DOA) estimation, referred to as two-step
knowledge-aided iterative estimation of signal parameters via
rotational invariance techniques (ESPRIT) method (Two-Step
KAI-ESPRIT), which achieves more accurate estimates than those of
prior art. The proposed Two-Step KAI-ESPRIT improves the estimation
of the covariance matrix of the input data by incorporating prior
knowledge of signals and by exploiting knowledge of the structure of
the covariance matrix and its perturbation terms. Simulation results
illustrate the improvement achieved by the proposed method.
\end{abstract}

%\begin{keywords}
%Large sensor arrays, knowledge-aided techniques, direction finding,
%high-resolution parameter estimation.
%\end{keywords}

\section{Introduction}

In array signal processing, direction-of-arrival (DOA) estimation is
a key task in a broad range of important applications including
radar and sonar systems, wireless communications and seismology
\cite{Vantrees,locsme,elnashar,manikas,cgbf,r19,scharf,bar-ness,pados99,
reed98,hua,goldstein,santos,qian,delamarespl07,delamaretsp,xutsa,xu&liu,
kwak,delamareccm,delamareelb,wcccm,delamarecl,delamaresp,delamaretvt,delamaretvt10,delamaretvt2011ST,
delamare_ccmmswf,jidf_echo,jidf,barc,lei09,delamare10,fa10,ccmavf,lei10,jio_ccm,
ccmavf,stap_jio,zhaocheng,zhaocheng2,arh_eusipco,arh_taes,rdrab,dcg,dce,dta_ls,
song,wljio,barc,saalt,mmimo,wence}. Classical high-resolution
methods for DOA estimation such as the multiple signal
classification (MUSIC) method \cite{schimdt}, the root-MUSIC
algorithm \cite{Barabell}, the estimation of signal parameters via
rotational invariance techniques (ESPRIT) \cite{Roy} and other
recent subspace techniques \cite{Steinwandt,Wang,Qiu} are based on
estimating the signal and noise subspaces from the sample covariance
matrix. The accuracy of the estimates of the covariance matrix is of
fundamental importance in parameter estimation. In practical
scenarios, only a limited number of samples is available and when
the number of samples is comparable to the number of sensor array
elements, the estimated and the true subspaces can significantly
diverge. This problem has been dealt with using random matrix theory
in \cite{Mestre1,Mestre2,Loubaton}, and the development of G-MUSIC,
which considers the asymptotic situation when both the sample size
and the number of array elements tend to infinity at the same rate.
It is then deduced that the introduced method more accurately
describes the case in which these two quantities are finite and
similar in magnitude \cite{Loubaton}.

In this work, we take into account a different approach to improve
the quality of the sample covariance matrix estimate in the finite
sample size region. Inspired by the structural approach of
\cite{Vorobyov1, Vorobyov2} and the use of prior knowledge about
signals \cite{Steinwandt2,Bouleux}, we develop an ESPRIT-based
technique that exploits both prior knowledge about signals and the
structure of the covariance matrix to improve the estimation
accuracy. Our approach determines the value of a scaling factor that
reduces the undesirable terms causing perturbations in the
estimation of the signal and noise subspaces in an iterative manner,
resulting in better estimates. This is done by choosing the set of
DOA estimates that have higher likelihood of being the set of true
DOAs. Furthermore, whereas in \cite{Vorobyov1,Vorobyov2}, the
undesirable terms are calculated based on steering vectors of
estimates, in the proposed method this computation makes use not
only of those estimates but also of the available prior knowledge of
DOAs. Considering a practical scenario, this task can be achieved
using signals coming from base stations or from static users in the
system.

The remainder of this paper is organized as follows. Section II describes the system model and its parameters. Section III presents the proposed Two-Step KAI-ESPRIT algorithm. In section IV, we deal with the computational complexity analysis by means of counting the multiplications and additions involved in the processing. Section V illustrates and discusses the simulation results and finally, the concluding remarks are drawn in Section VI.

\section{System Model and Background}
\label{sysmodel}

Let us assume that \textit{P} uncorrelated narrowband signals from
far-field sources impinge on a uniform linear array (ULA) of $M\ (M
> \textit{P})$ sensor elements from  directions ${\bm
\theta}=[\theta_{1},\theta_{2},\ldots, \theta_P]^T$. We also
consider that the sensors are spaced from each other by a distance $
d\leq\frac{\lambda_{c}}{2}$, where $\lambda_{c}$ is the signal wavelength,
and that without loss of generality, we have
\text{${\frac{-\pi}{2}\leq\theta_{1}\leq\theta_{2}\ldots
\leq\theta_P\leq \frac{-\pi}{2}}$}.

The $i$th data snapshot of the $M$-dimensional array output vector
can be modeled as
\begin{equation}
\bm x(i)=\bm A\,s(i)+\bm n(i),\qquad i=1,2,\ldots,N,
\label{model}
\end{equation}
where $\bm s(i)=[s_{1}(i),\ldots,s_{P}(i)]^T
\in\mathbb{C}^{\mathit{P\times1}}$ represents the zero-mean source
data vector, $\bm n(i) \in\mathbb{C}^{\mathit{M \times 1}}$ is the
vector of white circular complex Gaussian noise with zero mean and
variance $\sigma_n^2$, and $N$ denotes the number of available
snapshots. The Vandermonde matrix $\bm A(\bm \Theta)=[\bm
a(\theta_{1}),\ldots,\bm a(\theta_{P})] \in\mathbb
{C}^{\mathit{M\times P}}$, known as the array manifold, contains the
array steering vectors $\bm a(\theta_j)$ corresponding to the $n$th
source, which can be expressed as
\begin{equation}
\bm a(\theta_n)=[1,e^{j2\pi\frac{d}{\lambda_{c}}
\sin\theta_n},\ldots,e^{j2\pi(M-1)\frac{d}{\lambda_{c}}\sin\theta_n}]^T,
\label{steer}
\end{equation}
where $n=1,\ldots, P$. Using the fact that $\bm s(i)$ and $\bm n(i)$
are modeled as uncorrelated linearly independent variables, the
$M\times M$ signal covariance matrix is calculated by

\begin{equation}
\bm R=\mathbb E\left[\bm x(i)\bm x^H(i)
\right]=\bm A\,\bm R_{ss}\bm A^H+
\sigma_n^2\bm I_M,
\label{covariance}
\end{equation}
where the superscript \textit{H} and $\mathbb E[\cdot]$ in $\bm
R_{ss}=\mathbb E[\bm s(i)\bm s^H(i)]$ and in $\mathbb E[\bm n(i)\bm
n^H(i)]=\sigma_n^2\bm I_M^{}$ denote the Hermitian transposition and
the expectation operator and $\bm I_M$ stands for the $M\times M$
identity matrix. Since the true signal covariance matrix is unknown,
it must be estimated and a widely-adopted approach is the sample
average formula given by
\begin{equation}
 \bm {\hat{R}}=\frac{1}{N} \sum\limits^{N}_{i=1}\bm x(i)\bm x^H(i),
\label{covsample}
\end{equation}
whose estimation accuracy is dependent on the number of snapshots.

\section{Proposed Two-Step KAI-ESPRIT Algorithm}

In this section, we present the proposed two-step KAI-ESPRIT
algorithm and detail its main features.

We can start by expanding \eqref{covsample} using \eqref{model} as
follows:
\begin{eqnarray}
\bm {\hat{R}}=\frac{1}{N} \sum\limits^{N}_{i=1}(\bm A\,s(i)+\bm n(i))\:(\bm A\,s(i)+\bm n(i))^H \nonumber\\= \bm A\left\lbrace\frac{1}{N} \sum\limits^{N}_{i=1}\bm s(i)\bm s^H(i)\right\rbrace\bm A^H+\:\frac{1}{N} \sum\limits^{N}_{i=1}\bm n(i)\bm n^H(i)\;+\nonumber\\\bm A\left\lbrace\frac{1}{N} \sum\limits^{N}_{i=1}\bm s(i)\bm n^H(i)\right\rbrace\: +\:\left\lbrace\frac{1}{N} \sum\limits^{N}_{i=1}\bm n(i)\bm s^H(i)\right\rbrace\bm{A}^{H}
\label{expandedcovsample}
\end{eqnarray}
The first two terms of \text{$\bm {\hat{R}}$} in
\eqref{expandedcovsample} can be considered as estimates of the two
summands of \text{$\bm R$} given in  \eqref{covariance}, which
represent the signal and the noise components, respectively.  The
last two terms in \eqref{expandedcovsample} are undesirable
by-products, which  can be seen as estimates for the correlation
between the signal and the noise vectors. The system model under
study is based on noise vectors which are zero-mean  and also
independent of the signal vectors. Thus, the signal and noise
components are uncorrelated to each other. As a consequence, for a
large enough number of samples N, the last two  terms of
\eqref{expandedcovsample} tend to zero. Nevertheless, in practice
the number of available samples can be limited. In such situations,
the last two terms in \eqref{expandedcovsample} may have significant
values, which causes the deviation of the estimates of the signal
and the noise subspaces from the true signal and noise ones. The key
point of the proposed Two-Step KAI-ESPRIT algorithm is to modify the
sample data covariance matrix in the second step based on the
estimates obtained at the first step and the available known DOAs.
The modified covariance matrix is computed by deriving a scaled
version of the undesirable terms from $\bm {\hat{R}}$.

The steps of the proposed algorithm are listed in Table
\ref{Proposed_Two_Step_KA_ESPRIT}. The algorithm starts by computing
the sample data covariance matrix \eqref{covsample}. Next, the DOAs
are estimated using the ESPRIT algorithm. The superscript
$(\cdot)^{(1)}$ refers to the estimation task performed in the first
step. In the second step, the Vandermonde matrix is formed using the
 DOA estimates. Then, the amplitudes of
the sources are estimated such that the square norm of the
differences between the observation vector and the vector containing
estimates and the available known DOAs is minimized. This problem
can be formulated as
\begin{eqnarray}
    \hat{\bm{s}}(i)=\arg\min_{\substack{\bm
    s}}\parallel\bm{x}(i)-\hat{\bm{A}}\mathbf{s}\parallel^2_2.
    \label{minimization1}
\end{eqnarray}
The minimization of \eqref{minimization1} is achieved using the
least squares technique and the solution is described by
\begin{equation}
\hat{\bm{s}}(i)=(\mathbf{\hat{A}}^{H}\:\mathbf{\hat{A}})^{-1}\:\mathbf{\hat{A}}\:\bm{x}(i)
\label{minimization2}
\end{equation}
The noise component is then estimated as the difference between the
estimated signal and the observations made by the array, as given by
\begin{eqnarray}
 \hat{\bm n}(i)=\bm x(i)\:-\: \hat{\bm A}\:\hat{\bm s}(i).
\label{noise_component}
\end{eqnarray}
After estimating  the signal and noise vectors, the third term in
\eqref{expandedcovsample} can be computed as
\begin{align}
\bm{V}&\triangleq \hat{\bm{A}}\left\lbrace\frac{1}{N}
\sum\limits^{N}_{i=1}\bm \hat{\mathbf{s}}(i)\bm
\hat{\mathbf{n}}^H(i)\right\rbrace\nonumber\\&=\hat{\bm{A}}\left\lbrace\frac{1}{N}
\sum\limits^{N}_{i=1}(\mathbf{\hat{A}}^{H}\:\mathbf{\hat{A}})^{-1}\mathbf{\hat{A}}^{H}\bm{x}(i)\right.\nonumber\\&\left.\times(\bm{x}^{H}(i)-\bm{x}^{H}(i)\hat{\mathbf{A}}(\hat{\mathbf{A}}^{H}\hat{\mathbf{A}})^{-1}\:\hat{\mathbf{A}}^{H})\right\rbrace\nonumber\\&=\mathbf{\hat{Q}}_{A}\left\lbrace\frac{1}{N}
\sum\limits^{N}_{i=1}
\bm{x}(i)\bm{x}^H(i)\:\left(\mathbf{I}_{M}\:-\:\hat{\mathbf{Q}}_{A}\right)
\right\rbrace\nonumber\\&=\mathbf{\hat{Q}}_{A}\:\mathbf{\hat{R}}\:\mathbf{\hat{Q}}_{A}^{\perp},
\label{terms_deducted}
\end{align}
where
\begin{equation}
\mathbf{\hat{Q}}_{A}\triangleq \mathbf{\hat{A}}\:(\mathbf{\hat{A}}^{H}\:\mathbf{\hat{A}})^{-1}\:\mathbf{\hat{A}}^{H}
\end{equation}
is an estimate of the projection matrix of the signal subspace, and
\begin{equation}
\mathbf{\hat{Q}}_{A}^{\perp}\triangleq\mathbf{I}_{M}\:-\:\mathbf{\hat{Q}}_{A}
\end{equation}
is an estimate of the projection matrix of the noise subspace.

Lastly, the modified data covariance matrix is calculated by
computing a scaled version of the estimated terms from the initial
sample data covariance matrix as given
\begin{equation}
\label{modified_data_covariance}
\mathbf{\hat{R}}^{(n+1)} = \mathbf{\hat{R}}\:-\:\mathrm{\mu}\:(\mathbf{V}^{(n)}\:+\:\mathbf{V}^{(n)H}),\hspace{3mm}\textit{n}\hspace{1mm}=\hspace{1mm}\text{1}
\end{equation}
The scaling or reliability factor \text{$\mu $} increases from 0 to
1 incrementally, resulting in modified data covariance matrices.
Each of them gives origin to new estimated DOAs denoted by the
superscript  $(\cdot)^{(2)}$ by using the ESPRIT algorithm, as briefly described ahead.
Assuming the rank P is known, the eigenvalue decomposition (EVD) of the modified data covariance matrix \eqref{modified_data_covariance}   yields

\begin{equation}
\begin{array}{ccc}
\mathbf{\hat{R}}^{(n+1)}
\end{array} =\left[ \begin{array}{ccc}
\mathbf{\hat{U}}_{s} & \mathbf{\hat{U}}_{n}  \\
\end{array} \right]\left[ \begin{array}{ccc}
\mathbf{\hat{\Lambda}}_{s} & 0  \\
0 & \mathbf{\hat{\Lambda}}_{n}
\end{array} \right]\left[ \begin{array}{ccc}
\mathbf{\hat{U}}_{s}^{H} \\
\mathbf{\hat{U}}_{n}^{H}
\end{array} \right]
\label{Detalha_ESPRIT}
\end{equation}
where  the matrix $\mathbf{\hat{U}}_{s}$ $\in$ $\mathbb{C}^{\mathit{M\times P}}$ represents the signal subspace  and the matrix $\mathbf{\hat{U}}_{n}$ $\in$ $\mathbb{C}^{\mathit{M\times(M-P) }}$ represents the  the noise subspace respectively. The diagonal matrices $\mathbf{\hat{\Lambda}_{s}}$ and $\mathbf{\hat{\Lambda}_{n}}$ contain the \textit{P} largest and the \textit{M-P} smallest eigenvalues, respectively. We can form a twofold subarray configuration, as each row of the (Vandermonde) array steering matrix $\bm A(\bm \Theta)$ corresponds to one particular sensor element of the antenna array. The subarrays are specified by two $\mathit{(s\times M)}$-dimensional selection matrices $\mathbf{J_{1}}$  and $\mathbf{J_{2}}$ which choose $\mathit{s}$ elements of the $\mathit{M}$ existing sensors respectively, where $\mathit{s}$ is in the range  $\mathit{P\leq s <  M}$. For maximum overlap, the matrix $\mathbf{J_{1}}$ selects  the first $\mathit{s=M-1}$ elements and the matrix $\mathbf{J_{2}}$ selects the last $\mathit{s=M-1}$ rows of $\bm A(\bm \Theta)$.

Since the matrices $\mathbf{J_{1}}$  and $\mathbf{J_{2}}$ have now
been  computed, we can estimate the  operator $ \mathbf{\Psi} $ by
solving the approximation of the shift invariance equation
\eqref{shift_invariance_equation} given by
\begin{equation}
\mathbf{J}_{1}\:\mathbf{\hat{U}}_{s}\:\mathbf{\Psi}\:\approx\:\mathbf{J}_{2}\:\mathbf{\hat{U}}_{s}
\label{shift_invariance_equation}
\end{equation}
Using  the least square (LS) method, which yields
\begin{equation}
\hat{\mathbf{\Psi}}=\arg\min_{\substack{\mathbf
{\Psi}}}\parallel\mathbf{J}_{2}\:\mathbf{\hat{U}}_{s}\:-\:\mathbf{J}_{1}\:\mathbf{\hat{U}}_{s}\:\mathbf{\Psi}\parallel_{F}\:=\:\left(
\mathbf{J}_{1}\:\hat{\mathbf{U}}_{s}\right)
^{\dagger}\:\mathbf{J}_{2}\:\hat{\mathbf{U}}_{s},
\end{equation}
where $\parallel\cdot\parallel_{F}$ denotes the Frobenius norm and
$\left( \cdot\right)^{\dagger}$ stands for the pseudo-inverse.

Lastly, the eigenvalues $\lambda_{n}$ of $\hat{\mathbf{\Psi}}$
contain the estimates of the spatial frequencies $\gamma_{n}$
computed as
\begin{equation}
\gamma_{n}\:=\:\arg\left(\lambda_{n} \right)
\label{spatial_frequencies}
\end{equation}
so that the DOAs can be calculated as

\begin{equation}
\hat{\theta}_{n}\:=\:\arcsin\left(\frac{\gamma_{n}\:\lambda_{c}}{2\pi\:\mathit{d}} \right)
\label{doas_ESPRIT}
\end{equation}
where for \eqref{spatial_frequencies} and \eqref{doas_ESPRIT}
$\mathrm{n=1,\cdots,P}$.

Then, a new Vandermonde matrix is formed by the steering vectors of
those newly estimated DOAs and the available known DOAs. By using
the new matrix, it is possible to compute the new estimates of the
projection matrices of the signal \text{$ \mathbf{\hat{Q}}_{A}^{(n)}
$} and the noise \text{$ \mathbf{\hat{Q}}_{A}^{(n)\perp} $}
subspaces, both for n=2.

Next, employing the new estimates of the projection matrices, the
initial sample data matrix, $\bm {\hat{R}}$, and the number of
sensors and sources, the stochastic maximum likelihood objective
function \cite{Stoica} is computed for each value of \text{$\mu $},
as follows:
\begin{equation}
\mathit{U(\mu)}\,=\,\mathrm{ln\:det}\:\left(\mathbf{\hat{Q}}_{A}^{(2)}\:\mathbf{\hat{R}}\:\mathbf{\hat{Q}}_{A}^{(2)}\:+\:\dfrac{{\rm Trace}\{\mathbf{\hat{Q}}_{A}^{\perp\:(2)}\:\mathbf{\hat{R}}\}} {\mathrm{M-P}}\:\mathbf{\hat{Q}}_{A}^{\perp\:(2)}\right)
\label{SML_objective_function}
\end{equation}
The previous computation selects the set of unavailable DOA
estimates that have a higher likelihood. Then, the set of estimated DOAs corresponding to
the optimum value of \text{$\mu $} that minimizes
\eqref{SML_objective_function} is determined. Finally, the output of
the proposed Two-Step KAI-ESPRIT algorithm is formed by the set of
available known DOAs and the estimates of the unavailable DOAs, as
described in Table \ref{Proposed_Two_Step_KA_ESPRIT}.

%%%%%%%%%%%%%%%%%%%%%%%%%%%%%%%%%%%%%%%%%%%%%%%%%%%%%%%%%%%%%%%%%%%%%%%%%%
%
\begin{table}[htb!]
    \small
    \caption{\normalsize Proposed Two-Step KAI-ESPRIT Algorithm}
    \scalebox{1.0}\medskip{
        \begin{tabular}{|r l|}
            \hline

            \multicolumn{2}{|l|}{\small $\textbf{\underline{Inputs}:}$}\\[0.7ex]
            \multicolumn{2}{|l|}{\small$\mathit{M}$,\hspace{2mm}$\mathit{d}$,\hspace{2mm}$\lambda$,\hspace{2mm}$\mathit{N}$,\hspace{2mm}$\mathit{P}$ }\\[0.6ex]
            \multicolumn{2}{|l|}{\small\text{Received vectors}  $\bm x(1)$,\hspace{2mm}$\bm x(2)$,$\cdots$, $\bm x(N)$}\\[0.6ex]

            \multicolumn{2}{|l|}{\small \text{Prior knowledge $\rightarrow$ known DOAs:}\hspace{2mm}$\mathit{\theta}_{1}$,\hspace{2mm}$\mathit{\theta}_{2}$,
                $\cdots$,\hspace{2mm}$\mathit{\theta}_{q}$\hspace{2mm} $\mathit{1\leq q<P}$}\\[0.7ex]

            \multicolumn{2}{|l|}{\small $\textbf{\underline{Outputs}:}$}\\[0.6ex]

            \multicolumn{2}{|l|}{\small\text{Estimates}\hspace{1mm}$\mathit{\hat{\theta}_{q+1}^{(2)}}$,\hspace{2mm}$\mathit{\hat{\theta}_{q+2}^{(2)}}$,$\cdots$,\hspace{2mm}$\mathit{\hat{\theta}_{P}^{(2)}}$} \\[3.1ex]
            \hline

            \multicolumn{2}{|l|}{\small $\textbf{\underline{First step}:}$}\\[0.9ex]

            \multicolumn{2}{|l|}{\small $\mathbf{\hat{R}}=\frac{1}{N} \sum\limits^{N}_{i=1}\bm x(i)\bm x^H(i)$}\\[1.7ex]

            \multicolumn{2}{|l|}{\small $\{\mathit{\hat{\theta}_{1}}^{(1)},\:\mathit{\hat{\theta}_{2}}^{(1)},\cdots,\mathit{\hat{\theta}_{P}}^{(1)}\}\;\;\underleftarrow{ESPRIT}$ $\:(\mathbf{\hat{R}},P,d,\lambda)$}\\[1.1ex]

            \multicolumn{2}{|l|}{\small $\textbf{\underline{Second step}:}$}\\[0.9ex]

            \multicolumn{2}{|l|}{\small$\mathbf{\hat{A}}^{(1)}=\left[\mathbf{a}(\mathit{\hat{\theta}}_{1}^{(1)}),\mathbf{a}({\mathit{\hat\theta}_{2}^{(1)}}),\:\cdots,\mathbf{a}(\mathit{\hat{\theta}_{P}^{(1)}})\right]$} \\ [1.1ex]

            \multicolumn{2}{|l|}{\small \textbf{compute} for \textit{n}\hspace{1mm}=\hspace{1mm}\text{1}}\\[0.6ex]

            \multicolumn{2}{|l|}{\small $\mathbf{\hat{Q}}_{A}^{(n)}= \mathbf{\hat{A}}^{(n)}\:(\mathbf{\hat{A}}^{(n)H}\:\mathbf{\hat{A}}^{(n)})^{-1}\:\mathbf{\hat{A}}^{(n)H}$ \hspace{6mm}\text{(1)}}\\[0.9ex]

            \multicolumn{2}{|l|}{\small$\mathbf{\hat{Q}}_{A}^{(n)\perp}=\mathbf{I}_{M}\:-\:\mathbf{\hat{Q}}_{A}^{(n)}$ \hspace{27mm}\text{(2)}}\\[1.0ex]

            \multicolumn{2}{|l|}{\small $\mathbf{V}^{(n)}=\mathbf{\hat{Q}}_{A}^{(n)}\:\mathbf{\hat{R}}\:\mathbf{\hat{Q}}_{A}^{(n)\perp}$}\\[1.8ex]

            \multicolumn{2}{|l|}{\small \textbf{for} $\mathbf{\mu=}\hspace{1mm}\text{0}:\text{increment\hspace{1mm}:\hspace{1mm}1} $}\\[0.6ex]

            \multicolumn{2}{|l|}{\small \textbf{compute}\hspace{3mm}$ \mathbf{\hat{R}}^{(n+1)} = \mathbf{\hat{R}}\:-\:\mathrm{\mu}\:(\mathbf{V}^{(n)}\:+\:\mathbf{V}^{(n)H})$,\hspace{3mm}\textit{n}\hspace{1mm}=\hspace{1mm}\text{1}} \\[0.9ex]

            \multicolumn{2}{|l|}{\small $\{\mathit{\hat{\theta}_{q+1}}^{(2)},\:\mathit{\hat{\theta}_{q+2}}^{(2)},\cdots,\mathit{\hat{\theta}_{P}}^{(2)}\}\;\;\underleftarrow{ESPRIT}$ $\:(\mathbf{\hat{R}}^{(2)},\:P,d,\lambda)$}\\[1.1ex]

            \multicolumn{2}{|l|}{\small$\mathbf{\hat{A}}^{(2)}=\left[\mathbf{a}(\mathit{\theta}_{1}), \ldots, \mathbf{a}(\mathit{\theta}_{q}),\mathbf{a}(\mathit{\hat{\theta}}_{q+1}^{(2)}),\ldots,\mathbf{a}(\mathit{\hat{\theta}_{P}^{(2)}})\right]$} \\[1.4ex]

            \multicolumn{2}{|l|}{\small \textbf{repeat}\hspace{1mm}\text{(1) and (2),}\hspace{1mm} \textit{n}\hspace{1mm}=\hspace{1mm}\text{2}}\\[0.6ex]

            \multicolumn{2}{|l|}{\small $\mathit{U(\mu)}\,=\,\mathrm{ln\:det}\:\left(\mathbf{\hat{Q}}_{A}^{(2)}\:\mathbf{\hat{R}}\:\mathbf{\hat{Q}}_{A}^{(2)}\:+\:\dfrac{{\rm Trace}\{\mathbf{\hat{Q}}_{A}^{\perp\:(2)}\:\mathbf{\hat{R}}\}} {\mathrm{M-P}}\:\mathbf{\hat{Q}}_{A}^{\perp\:(2)}\right),$}\\[1.1ex]

            \multicolumn{2}{|l|}{\small \textbf{end for}}\\[1.0ex]

            \multicolumn{2}{|l|}{\small $\mathit{\mu}_{\rm opt}=\arg \min \hspace{1mm}\mathit{U(\mu)}$}\\[1.0ex]

            \multicolumn{2}{|l|}{\small $\mathrm{DOAs}=  \{\mathit{\theta}_{1},\cdots,\:\mathit{\theta}_{q},\mathit{\hat{\theta}_{q+1}({\mu}_{opt})},\cdots,\mathit{\hat{\theta}_{P}({\mu}_{opt})}\}$}\\[1.0ex]

            \hline
        \end{tabular}
    }
    %\caption{Table to test captions and labels}
    \label{Proposed_Two_Step_KA_ESPRIT}

\end{table}

\section{Computational Complexity Analysis}
\label{Comput_Complexity}

In this section, we evaluate the computational cost of the proposed
Two-Step KAI-ESPRIT algorithm which is compared to the following
subspace methods: MUSIC, Root-MUSIC, Conjugate Gradient (CG) and
Auxiliary Vector Filtering (AVF). All of them made use of Singular
Value Decomposition (SVD) of the sample covariance matrix
\eqref{covsample}. The computational cost in terms of number of
multiplications and additions is depicted in Tables
\ref{Comput_Complexity1} and \ref{Comput_Complexity2}, where
$\mathrm{\Delta} $ is the search step and
$\mathrm{\tau}=\frac{1}{increment} +1$, and this increment is
defined in Table \ref{Proposed_Two_Step_KA_ESPRIT}.

As can be seen, the proposed algorithm shows a relatively high
computational burden with $\mathit{O(\tau(3M^{3}))}$, where  $\tau$
is typically an integer inside [1 20]. This is motivated  by two
reasons. The first is  that the modified data covariance matrix
\eqref{modified_data_covariance} needs to be computed  $\tau$ times.
The second is the need for three matrix multiplications of order
$\left[ M\times M\right] $ that define  the undesirable
by-products\eqref{expandedcovsample}, \eqref{terms_deducted} to be
subtracted from the sample covariance matrix
\eqref{modified_data_covariance}. A brief comparison  between  the
computational cost of the proposed algorithm and the others listed
in Table \ref{Comput_Complexity1} can be done considering the
dominant terms in the required multiplications. For this purpose we
suppose typical values of the increment  $\tau=20$ and of the search
step $\Delta=0.1$ degrees.  For MUSIC, CG and AVF, which require
peak searches, the  comparisons yield

 \begin{equation}
  \frac{\mathcal{O}\mathrm{\left(Two-Step\: KAI\: ESPRIT\right)}} {\mathcal{O}\mathrm{\left(MUSIC\right)}}\approx\frac{60\:M^{3}\Delta}{180\:M^{2}}\:\approx\:\frac{M}{30}
  \label{order_MUSIC}
 \end{equation}

 \begin{equation}
 \frac{\mathcal{O}\mathrm{\left(Two-Step\: KAI\: ESPRIT\right)}} {\mathcal{O}\mathrm{\left(CG\right)}}\approx\frac{60\:M^{3}\Delta}{180\:M^{2}\:P}\approx\:\frac{M}{30\:P}
 \label{order_CG}
 \end{equation}

 \begin{equation}
 \frac{\mathcal{O}\mathrm{\left(Two-Step\: KAI\: ESPRIT\right)}} {\mathcal{O}\mathrm{\left(AVF\right)}}\approx\frac{60\:M^{3}\Delta}{180\:M^{2}\:3P}\approx\:\frac{M}{90\:P}
 \label{order_AVF}
 \end{equation}

For root-MUSIC and the original ESPRIT, which do not involve peak
searches, we have
 \begin{equation}
 \frac{\mathcal{O}\mathrm{\left(Two-Step\: KAI\: ESPRIT\right)}} {\mathcal{O}\mathrm{\left(root-MUSIC\right)}}\approx \frac{60\:M^{3}}{2\:M^{3}}\approx\:30
 \label{order_root_MUSIC}
 \end{equation}

 \begin{equation}
 \frac{\mathcal{O}\mathrm{\left(Two-Step\: KAI\: ESPRIT\right)}} {\mathcal{O}\mathrm{\left(original\:ESPRIT \right)}}\approx \frac{60\:M^{3}}{2\:M^{2}\:P}\approx\:\frac{30M}{P}
 \label{order_original_ESPRIT}
 \end{equation}
 \medskip

As can be noticed in
\eqref{order_MUSIC},\eqref{order_CG},\eqref{order_AVF},\eqref{order_root_MUSIC}
and \eqref{order_original_ESPRIT}, the order of multiplications of
the proposed algorithm  applied to P=4 supposed signals impinging a
ULA formed with M=40 sensors is approximately $1.3\times$ the order
of multiplications of the MUSIC, $0.3\times$ the order of
multiplications of the CG and $0.1\times$ the order of
multiplications of the AVF algorithms, respectively. Therefore, in
this particular scenario, the number of multiplications of Two-Step
KAI ESPRIT can be considered to be approximately equal or less the
number of multiplications of the algorithms that require peak search
which are considered in this work. It can also be seen that the
number of multiplications of the proposed algorithm is roughly
$30\frac{M}{P} \times$ the number of the multiplications of the
original ESPRIT and $(30\frac{M}{P} \times$ the number of
multiplications of root-MUSIC, which do not require extensive peak
searches. To sum up, the  comparison with the algorithms that
require peak search allow us to consider that the relatively high
computational burden, which is associated with extra matrix
multiplications and the increment $\tau$ applied to cancel the
undesirable by-products, is not a too high cost to be paid for the
improved performance achieved. Similar results can be obtained for
the order of additions, except for the comparison with the AVF
algorithm, which yields

   \begin{equation}
   \frac{\mathcal{O}\mathrm{\left(Two-Step\: KAI\: ESPRIT\right)}} {\mathcal{O}\mathrm{\left(AVF\right)}}\approx\frac{60\:M^{3}\Delta}{180\:M^{2}\:4P}\approx\:\frac{M}{90\:P}
   \label{order_addition_AVF}
   \end{equation}
   what means that in  the same scenario described before the order of multiplications of our proposed algorithm
   is approximately $(0.08\times)$ the order of multiplications of the AVF algorithm.

\begin{table}[!h]
    \caption{Computational complexity applying the SVD}
    \vspace{2mm}
    \centering
    \begin{tabular}{|l|l|p{6cm}|}
        \hline
        Algorithm & Multiplications \\[2pt]
        \hline
        MUSIC \cite{schimdt} & $\mathrm{\frac{180}{\Delta}[M^{2}+M(2-P)-P]+8MN^{2}}$  \\[2pt]
        \hline
        root-MUSIC\cite{Barabell} & $\mathrm{2M^{3}-M^{2}P+8MN^{2}}$  \\[2pt]
        \hline
        AVF \cite{Grover} & $\mathrm{\frac{180}{\Delta}[M^{2}(3P+1)+M(4P-2)+P+2]}$\\[1pt]
        &$\mathrm{+M^{2}N}$  \\[2pt]
        \hline
        CG \cite{Semira} & $\mathrm{\frac{180}{\Delta}[M^{2}(P+1)+M(6P+2)+P+1]+M^{2}N}$  \\[2pt]
        \hline
        ESPRIT\cite{Roy} & $\mathrm{2M^{2}P+M(P^{2}-2P+8N^{2})+8P^{3}-P^{2} }$  \\[2pt]
        \hline
        Two-Step & $\mathrm{\tau[3M^{3}+M^{2}(3P+2)+M(\frac{5}{2}P^{2}-\frac{3}{2}P+8N^{2})}$  \\[1pt]

        KAI-ESPRIT&$\mathrm{+P^{2}(\frac{17}{2}P+\frac{1}{2})+1]}$ \\[1pt]

        (Proposed)& $\mathrm{+[2M^{3}+M^{2}(3P)+M(\frac{5}{2}P^{2}-\frac{3}{2}P+8N^{2})}$ \\[1pt]

        & $\mathrm{+P^{2}(\frac{17}{2}P+\frac{1}{2})]}$\\[2pt]

        \hline
    \end{tabular}

    \label{Comput_Complexity1}
\end{table}
\begin{table}[!h]
    \caption{Computational complexity applying the SVD}
    \vspace{2mm}
    \centering
    \begin{tabular}{|l|l|p{6cm}|}
        \hline
        Algorithm & Additions \\[2pt]
        \hline
        MUSIC \cite{schimdt} & $\mathrm{\frac{180}{\Delta}[M^{2}+M(1-P)-2]+8MN^{2}}$  \\[2pt]
        \hline
        root-MUSIC\cite{Barabell} & $\mathrm{2M^{3}-M^{2}P+M(8N^{2}-2)+1}$  \\[2pt]
        \hline
        AVF \cite{Grover} & $\mathrm{\frac{180}{\Delta}[4PM^{2}+M(2P-3)-3P+2]}$\\[1pt]
        &$\mathrm{+ M^{2}(N-1)}$  \\[2pt]
        \hline
        CG \cite{Semira} & $\mathrm{\frac{180}{\Delta}[M^{2}(P+1)+M(5P+1)-3P-2]}$\\[1pt]
        &$\mathrm{+M^{2}(N-1)}$  \\[2pt]
        \hline
        ESPRIT\cite{Roy} & $\mathrm{2M^{2}P+M(P^{2}-4P+8N^{2})+8P^{3}}$\\[1pt]
        &$\mathrm{-P^{2}+2P }$  \\[2pt]
        \hline
        Two-Step & $\mathrm{\tau[3M^{3}+M^{2}(3P-1)}$\\[1pt]

        &$\mathrm{+M(\frac{5}{2}P^{2}-\frac{9}{2}P+8N^{2}+2)}$  \\[1pt]

        KAI-ESPRIT&$\mathrm{+P(8P^{2}-2P-\frac{5}{2})]}$\\[1pt]

        (Proposed)& $\mathrm{+[2M^{3}+M^{2}(3P-3)}$\\[1pt]
        &$\mathrm{+M(\frac{5}{2}P^{2}-\frac{9}{2}P+8N^{2}+1)}$ \\[1pt]

        & $\mathrm{+P(8P^{2}-2P-\frac{5}{2})]}$\\[2pt]

        \hline
    \end{tabular}

    \label{Comput_Complexity2}
\end{table}

\section{Simulations}
\label{simulations}

\begin{figure}[!h]

    \centering % para centralizarmos a figura
    \includegraphics[width=9cm,scale=1.0]{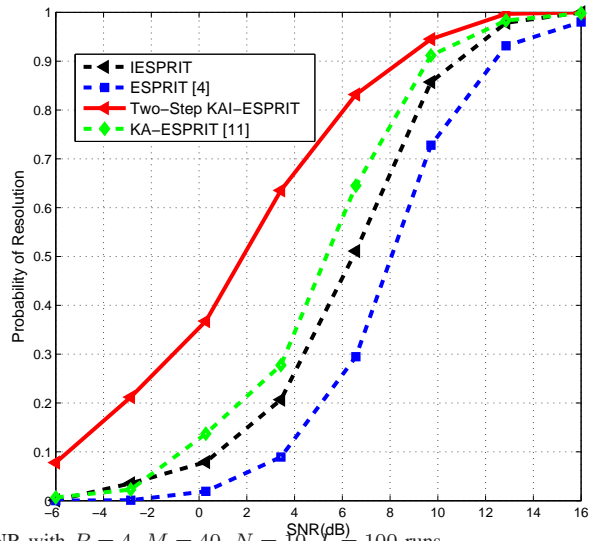} % leia abaixo
    \vspace{-2.5em}\caption{Probability of resolution versus SNR with $P=4$, $M=40$, $N=10$, $L=100$ runs}
    \label{figura:WSA_PR_2deg_40sens_10snap_100runs}
\end{figure}

In this section, we examine the performance  of the proposed
Two-Step KAI-ESPRIT algorithm in terms of probability of
resolution and RMSE and compare them to the conventional ESPRIT
\cite{Roy}, the Iterative ESPRIT (IESPRIT), which is also developed
here by combining the approach in \cite{Vorobyov2} that exploits
knowledge of the structure of the covariance matrix and its
perturbation terms and the standard ESPRIT, and the Knowledge-Aided
ESPRIT (KA-ESPRIT) \cite{Steinwandt2} using general linear
combination \cite{Stoica2}. We employ a ULA with \textit{ M=40}
sensors, inter-element spacing $\Delta=\frac{\lambda_{c}}{2}$ and
assume there are four uncorrelated complex Gaussian signals with
equal power impinging on the array. The closely-spaced sources are
separated by \textit{$2^{o}$}, at $\mathrm
(13^{o},15^{o},17^{o},19^{o})$,  and the number of available
snapshots is \textit{N}=10. Furthermore, we presume a priori
knowledge of the two last true DOAs at $\mathrm (17^{o},19^{o})$
only in the proposed Two-Step KAI-ESPRIT and in the KA-ESPRIT. In
Fig. \ref{figura:WSA_PR_2deg_40sens_10snap_100runs}, we show the
probability of resolution versus SNR. We take into account the
criterion \cite{Stoica3}, in which two sources with DOA $\theta_{1}$
and $\theta_{2}$  are said to be resolved if their respective
estimates $\hat{\theta}_{1}$ and $\hat{\theta}_{2}$ are such that
both $\left|\hat{\theta}_{1} -\theta_{1}\right|$ and
$\left|\hat{\theta}_{2} -\theta_{2}\right|$ are less than
$\left|\theta_{1} -\theta_{2}\right|/2$. The proposed Two-Step
KAI-ESPRIT algorithm outperforms KA-ESPRIT \cite{Roy,Steinwandt2},
IESPRIT and the standard ESPRIT. %It must be noted that, like the
%original ESPRIT, its iterative variation IESPRIT is not
%knowledge-aided.

\begin{figure}[!h]

    \centering % para centralizarmos a figura
    \includegraphics[width=9cm, scale=1.0]{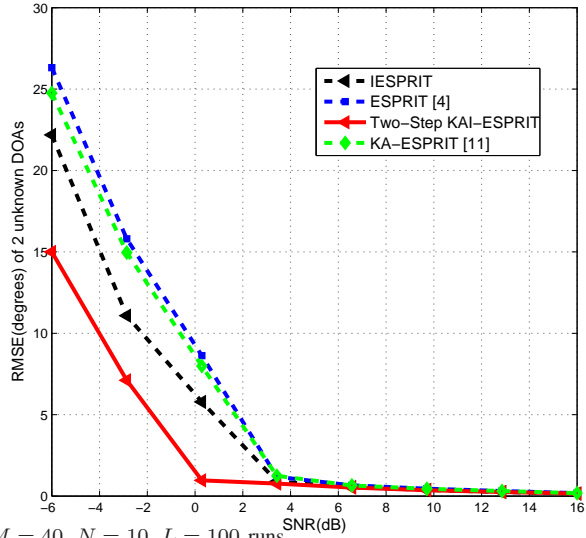} % leia abaixo
    \vspace{-2.5em}\caption{RMSE versus SNR with $P=4$, $M=40$, $N=10$, $L=100$ runs}
    \label{figura:WSA_RMSE_2deg_40sens_10snap_100runs}
\end{figure}

\begin{figure}[!h]

    \centering % para centralizarmos a figura
    \includegraphics[width=9cm, scale=1.0]{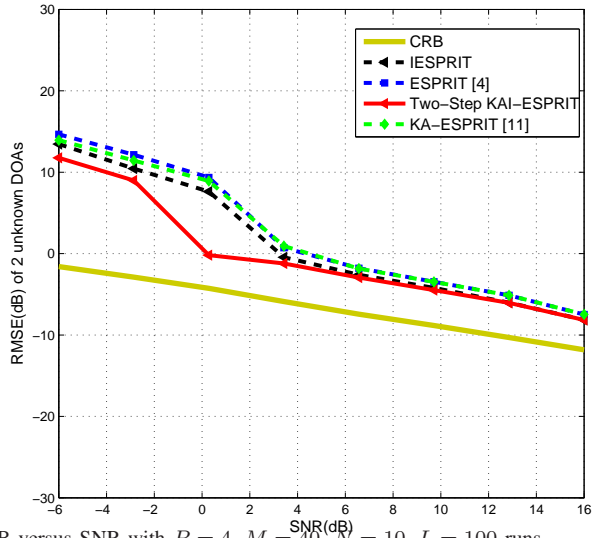} % leia abaixo
    \vspace{-2.5em}\caption{RMSE and the square root of CRB versus SNR with $P=4$, $M=40$, $N=10$, $L=100$ runs}
    \label{figura:WSA_RMSE_CRB_2deg_40sens_10snap_100runs}
\end{figure}

In Fig. \ref{figura:WSA_RMSE_2deg_40sens_10snap_100runs}, it is
shown the RMSE of the two supposedly unknown DOAs versus SNR. For
the computations we adopted the expression
\begin{equation}
\centering \mathrm{RMSE}
=\sqrt{\frac{1}{L\:P}\sum\limits^{L}_{l=1}\sum\limits^{P}_{p=1}\bm(\theta_p
-\bm \hat{\theta}_p(l))}, \label{RMSE_run}
\end{equation}
where L=number of trials. Alternatively, Fig.
\ref{figura:WSA_RMSE_2deg_40sens_10snap_100runs} can be expressed in
terms of dB as shown in Fig.
\ref{figura:WSA_RMSE_CRB_2deg_40sens_10snap_100runs}, where  the
term CRB refers to the square root of the deterministic Cramér-Rao
bound \cite{Stoica4}. The results show the superior performance of
the proposed Two-Step KAI-ESPRIT algorithm in terms of RMSE and also
of probability of resolution. In particular, the proposed technique
can obtain a higher probability of resolution and a lower RMSE than
existing techniques for lower values of SNR. According to Figs.
\ref{figura:WSA_PR_2deg_40sens_10snap_100runs} and
\ref{figura:WSA_RMSE_2deg_40sens_10snap_100runs}, Two-Step
KAI-ESPRIT can obtain for the same performance in RMSE or
probability of resolution gains in SNR that range from $0.5$ to
$3.0$ dB.

\section{Conclusions}

We have proposed in this work the Two-Step KAI-ESPRIT algorithm
which exploits prior knowledge of source signals and the structure
of the covariance matrix and its perturbations. The proposed
Two-Step KAI-ESPRIT algorithm can obtain significant gains in RMSE
or probability of resolution performance over previously reported
techniques, and has excellent potential for applications with short
data records in large-scale antenna systems for wireless
communications, radar and other large sensor arrays. The relatively
high computational burden required, which is associated with extra
matrix multiplications and the increment applied to reduce the
undesirable by-products can be justified for the superior
performance achieved. Future work will consider approaches to
reducing computational cost.

\end{document}